\documentclass[12pt,preprint,nofootinbib,nobibnotes,nonumbib]{article}
\usepackage{color}
\usepackage{amsmath} 
\usepackage{amsxtra} 
\usepackage{amscd}
\usepackage{amsthm} 
\usepackage{amsfonts} 
\usepackage{amssymb} 
\usepackage{eucal} 

\definecolor{Cyan}{rgb}{0,0,1}
\color{Cyan}

\begin{document}

\title{$C^\infty$ Functions on the Stone-${\check {\rm C}}$ech 
Compactification of the Integers}

\author{Larry B. Schweitzer}

\date{October 2014}

\maketitle

\begin{abstract}   
We construct an algebra $A=\ell^{\infty \infty}({\Bbb Z})$ of 
smooth functions which is dense 
in the pointwise multiplication algebra $\ell^\infty({\Bbb Z})$
of sup-norm bounded functions on the integers $\Bbb Z$.  
The algebra $A$ properly contains the sum of the algebra 
$A_c=\ell_c^\infty({\Bbb Z})$ and the ideal ${\cal S}({\Bbb Z})$, 
where $A_c$ is the algebra of finite linear 
combinations of projections in $\ell^\infty({\Bbb Z})$ 
and ${\cal S}({\Bbb Z})$ is the
pointwise multiplication algebra of Schwartz functions.  
The algebra $A$ is characterized as the set of functions whose 
\lq\lq first derivatives\rq\rq\ vanish 
rapidly at each point in the Stone-${\check {\rm C}}$ech 
compactification of $\Bbb Z$.

\noindent{
2000 Mathematics Subject Classification: 
%NCDG:
58B34,
%Rings & Algebras of continuous, differentiable or analytic functions:
46E25, 
%Special constructions of spaces (from ultrafilters, etc):
54D80, 
%Extensions of spaces (compactifications, etc):
54D35, 
%Extremely disconnected spaces, F-spaces, etc:
54G05, 
%Sequence spaces (including Kothe sequence spaces):
46A45.}
\end{abstract}

\tableofcontents

\vskip\baselineskip
\section{Introduction}

In a previous paper [Sch, 1998], a notion of smooth functions on 
the Cantor set was developed.  Recall that a {\it totally disconnected}
topological space has a basis of clopen sets.  The Cantor set is such a space.
In this paper, we attempt to construct smooth functions on the
Stone-${\check {\rm C}}$ech 
compactification of the integers.  In addition to being
totally disconnected, this space is {\it extremely disconnected}, 
which means that the closure of every open set is clopen. 

We will be working with the $C^\star$-algebra of all bounded complex-valued
functions (or sequences) on the integers ${\Bbb Z}$, under pointwise
multiplication, with pointwise complex-conjugation for involution.  
We denote this algebra by $\ell^\infty({\Bbb Z})$.  None of the theorems
will use the additive and multiplicative structure of ${\Bbb Z}$, so that
any countable discrete set can be substituted for ${\Bbb Z}$.  
However, when specific
objects are constructed in the examples, we may make reference to the
underlying set of integers.   Recall that the integers ${\Bbb Z}$ consists of 
all whole numbers from $-\infty$ to $\infty$, whereas the natural
numbers ${\Bbb N}$ contains only the whole numbers from $0$ to $\infty$.    
The norm on $\ell^\infty({\Bbb Z})$ 
is the sup-norm $\| \quad \|_\infty$, defined by
$\| \varphi \|_\infty = \sup_{n \in {\Bbb Z}} |\varphi(n)|$.

We will use the notation $c_0({\Bbb Z})$ for the complex-valued functions
(sequences) on ${\Bbb Z}$ which vanish at infinity.

\vskip\baselineskip
\section{The Stone-${\check {\rm C}}$ech Compactification of $\Bbb Z$}

We recall the standard definition of the 
Stone-${\check {\rm C}}$ech compactifiaction
from [Roy, 1968].
Let $\bf F$ be the set of all real-valued functions from $\Bbb Z$
into the closed interval $I=[-1,1]$.  (So $\bf F$ is the set of
real-valued functions in the unit ball of $\ell^\infty({\Bbb Z})$.)  
Let $${\bf X} = \prod_{\bf F} I$$
be the $\bf F$-fold cartesian product of unit intervals $I$.
By the Tychonoff theorem [Roy, 1968, Chapter 9, Theorem 19], 
this is a compact Hausdorff space.
Let 
	$$i:{\Bbb Z} \hookrightarrow {\bf X}  
	{\rm \quad where \quad}
	i(n)\longmapsto\bigl\{f(n)\bigr\}_{f \in \bf F}$$
be the natural inclusion map.   

\noindent{\bf Definition 1.1. \ }  
We define the {\it Stone-${\check {\rm C}}$ech 
Compactification of ${\Bbb Z}$}, denoted by $\beta({\Bbb Z})$, to be the closure 
of the image $i({\Bbb Z})$ in ${\bf X}$.   We let $C(\beta({\Bbb Z}))$ denote
the continuous complex-valued functions on $\beta({\Bbb Z})$.
See [Roy, 1968] for the basic properties of $\beta({\Bbb Z})$.

\noindent{\bf Definition 1.2. \ }
Let $n_0 \in \beta({\Bbb Z})$, and let $\{n_\alpha\}_{\alpha \in \Lambda}$
be a net converging to $n_0$, where $\Lambda$ is some
directed set, and $n_\alpha \in {\Bbb Z}$ for each $\alpha \in \Lambda$.
Let $\varphi \in \ell^\infty({\Bbb Z})$.  
Then {\it $\varphi(n_0)$} is defined as the limit
of the net $\{\varphi(n_\alpha)\}_{\alpha \in \Lambda}$.    
In this way $\varphi$ defines a continuous function
on $\beta({\Bbb Z})$, giving an isomorphism of commutative $C^\star$-algebras
$\ell^\infty({\Bbb Z}) \cong C(\beta({\Bbb Z}))$.   
(The fact that $\varphi$ is continuous,
and that the limit defining $\varphi(n_0)$ converges, can be verified
directly from the definition of 
the Stone-${\check {\rm C}}$ech compactification above,
using the fact that the real and imaginary parts of 
$\varphi/{\| \varphi \|_\infty}$ are functions in $\bf F$.)

Let $p\in \ell^\infty({\Bbb Z})$ be a projection.  
In other words for each $n\in {\Bbb Z}$, 
$p(n) = 0$ or $p(n)=1$.  Then $p$ also 
defines a projection on $\beta({\Bbb Z})$.
In fact, if $n_0\in \beta({\Bbb Z})$ and 
$\{n_\alpha\}_{\alpha \in \Lambda}$ are as above, then
$p(n_\alpha)$ is either eventually equal to $1$ (in the case $p(n_0)=1$) or 
eventually equal to $0$ (in the case $p(n_0)=0$).  
Thus each point $n_0
\in \beta({\Bbb Z})$ defines a family ${\cal F}_{n_0}$ of subsets of $\Bbb Z$, 
where $S \in {\cal F}_{n_0}$ if and only if the projection $p$ whose support is 
equal to $S$ satisfies $p(n_0)=1$.

\noindent{\bf Definition 1.3. \ }  
A family of subsets $\cal F$ of ${\Bbb Z}$
is a {\it filter} on ${\Bbb Z}$ if it is closed under finite intersections
$$S, T \in {\cal F} \Longrightarrow S\cap T \in {\cal F}\eqno{\rm (1.4a)} $$
and supersets
$$S \in {\cal F} \quad{\rm and}\quad
T \supset S
\Longrightarrow T \in {\cal F},\eqno{\rm (1.4b)} $$
and if the empty set $\emptyset$ in {\it not} in ${\cal F}$.
A filter $\cal U$ on ${\Bbb Z}$ is an 
{\it ultrafilter} on ${\Bbb Z}$ if
$$S \subseteq {\Bbb Z} \Longrightarrow S \in {\cal F} 
\quad{\rm or}\quad  S^{c} \in {\cal F}.\eqno{\rm (1.4c)} $$
The map $n_0
\in \beta({\Bbb Z}) \mapsto {\cal F}_{n_0}$ from the previous paragraph 
defines a map {\it ${\cal UF}$} from $\beta({\Bbb Z})$ to the set 
of ultrafilters on ${\Bbb Z}$.  (Use the definition of 
$\beta({\Bbb Z})$ to check this.)

%If $\{ n_\alpha \}_{\alpha \in \Lambda}$
%and $\{ m_\alpha \}_{\alpha \in \Sigma}$ are two nets in ${\Bbb Z}$ converging
%to the same point $n_0$ in $\beta({\Bbb Z})$,  they are eventually in the
%same subsets of ${\Bbb Z}$.  Thus the map $n_0 \mapsto {\cal UF}(n_0)$
%is a well-defined map $\beta({\Bbb Z}) \rightarrow {\cal UF}$.

A {\it principal filter} is of the form ${\cal U}_n = <S \subseteq {\Bbb Z} \mid
n \in S>$ for some $n \in {\Bbb Z}$.  Such a filter is also an ultrafilter.
The image of ${\Bbb Z} \subseteq \beta({\Bbb Z})$ under the map ${\cal UF}$
is precisely the set of principal ultrafilters on ${\Bbb Z}$.

We show that the map ${\cal UF}$ is an isomorphism by constructing an inverse
map ${\cal FU}$.  Let ${\cal U}$ be any ultrafilter on ${\Bbb Z}$.
Define a directed set $\Lambda$ to be ${\cal U}$ with superset order.
That is $\alpha \le \beta \Longleftrightarrow \beta \subseteq \alpha$.
Thus smaller sets are \lq\lq bigger\rq\rq\ in this order.  
In the case of the principal
ultrafilter $\Lambda={\cal U}_n$, the singleton $\{n\}$ is the biggest
element.  In general, $\Lambda$ has a biggest element if and only if 
it comes from a principal ultrafilter.
Next, for $\alpha\in \Lambda$ choose any $n_\alpha\in \alpha \subseteq {\Bbb Z}$.
We show that $\{ n_\alpha \}_{\alpha \in \Lambda}$ converges to an 
element of $\beta({\Bbb Z})$.

Let $\varphi \in {\bf F}$, where ${\bf F}$ is the defining family of 
functions for $\beta({\Bbb Z})$ from Definition 1.1.  
We wish to show that $\{ \varphi(n_\alpha) \}_{\alpha \in \Lambda}$
converges.  Define
\begin{equation}
\varphi_{+}(n) = 
\begin{cases}
\varphi(n) &\text{if $\varphi(n)\ge 0$;}\\
			0 &\text{otherwise.}
\end{cases}
\nonumber
\end{equation}
and
\begin{equation}
\varphi_{-}(n) = 
\begin{cases}
-\varphi(n) &\text{if $\varphi(n)\le 0$;}\\
			0 &\text{otherwise.}
\end{cases}
\nonumber
\end{equation}
Then $\varphi = \varphi_{+} - \varphi_{-}$, and $\varphi_{+},
\varphi_{-} \in {\bf F}$.   It suffices to show that 
$\{ \varphi_{+} (n_\alpha) \}_{\alpha \in \Lambda}$
and $\{ \varphi_{-} (n_\alpha) \}_{\alpha \in \Lambda}$
each converge.  So without loss of generality, we take $\varphi$
with range in the unit interval $[0,1]$.

If $\varphi$ were a projection $p$, we would be done.  Simply let $S$
be the support of $p$.  If $S \in  {\cal U} = \Lambda$, then 
$\{ p (n_\alpha) \}_{\alpha \in \Lambda}$ is eventually $1$
(when $\alpha \ge S$). Otherwise,
it is eventually $0$, and in either case they converge.  
We proceed by writing $\varphi$ as an infinite series of projections.
Define the set of integers
$$S_{1\over 2} = \{ n \in {\Bbb Z} \mid {1\over 2} \le \varphi(n) \le 1 \}.$$
Let $p_{1\over 2}$ be the projection corresponding to $S_{1\over 2}$.
Then $\varphi - {1\over 2}p_{1\over 2}$ has its range in the interval
$[0,{1\over 2}]$.  We repeat the process to get a new function with 
range in the interval $[0,{1\over 4}]$, etc, until we get:
$$\varphi = {1\over 2}p_{1\over 2} + {1\over 4}p_{1\over 4} + \dots 
	= \sum_{q=1}^{\infty} {1\over 2^q} p_q, \eqno(1.5)$$
an infinite series that converges absolutely (and geometrically fast) 
in sup norm to $\varphi$.  To see that 
$\{ \varphi (n_\alpha) \}_{\alpha \in \Lambda}$ converges,
now use a standard series argument.  Let $\epsilon>0$ be given, and
find $N$ large enough so that the sup norm of the tail 
$$ \biggl{\Vert} \sum_{q=N}^{\infty} {1\over 2^q} p_q \biggr{\Vert}_{\infty}$$
is less than $\epsilon$.  Then find $\alpha \in \Lambda$ sufficiently 
large so that
$\alpha \subseteq S_{1/ 2^q}$ or $\alpha \subseteq S^c_{1/ 2^q}$
for each $q=1,\dots N$.  (One could first find an $\alpha$ that works
for each $S_{1/ 2^q}$ separately, using the ultrafilter condition (1.4c). 
Then, find an $\alpha$ that works for all $S_{1/ 2^q}, q=1,\dots N$ 
using the finite intersection property (1.4a).)   For any $\beta \in \Lambda$
beyond this $\alpha$, the first $N$ projections in (1.5) have settled
down to their final value (either $0$ or $1$) on the net
$\{ n_\alpha \}_{\alpha \in \Lambda}$.  Thus  the net of real numbers
$\{ \varphi (n_\alpha) \}_{\alpha \in \Lambda}$ converges.
Since $\{ n_\alpha \}_{\alpha \in \Lambda}$ was an arbitrary net from 
an arbitrary ultrafilter on ${\Bbb Z}$,  this shows that ${\cal FU}$
is a well-defined map from the ultrafilters on ${\Bbb Z}$ into $\beta({\Bbb Z})$.
One easily checks that the compositions ${\cal UF} \circ {\cal FU}$ and 
${\cal FU} \circ {\cal UF}$ are identity maps, and so
we have proved:

\noindent{\bf Proposition 1.6. \ }  
The map ${\cal UF}$ is an isomorphism 
of the Stone-${\check {\rm C}}$ech 
compactification $\beta({\Bbb Z})$ with the set of ultrafilters on ${\Bbb Z}$.
Under this map, a point $n_0\in \beta({\Bbb Z})$ is taken to the ultrafilter
of sets that any net of integers converging to $n_0$ is eventually in.

\vskip\baselineskip
\section{Definition of the Smooth Functions $\ell^{\infty \infty}({\Bbb Z})$ 
( also denoted by $C^\infty(\beta({\Bbb Z}))$ )}

\noindent{\bf Definition 2.1. \ }  
Define {\it $\ell^\infty_c({\Bbb Z})$} to 
be the finite span of projections in $\ell^\infty({\Bbb Z})$.
This is a dense $\star$-subalgebra of $\ell^\infty({\Bbb Z})$, and plays an 
analogous role to $\ell^\infty({\Bbb Z})$ 
as $c_c(\Bbb Z)$, the compact
(or finite) support functions on $\Bbb Z$, does to $c_0({\Bbb Z})$.  
(The series expansion (1.5) proves the density.) 

For $n_0 \in \beta({\Bbb Z})$, choose a net $\{ n_\alpha \}_{\alpha \in \Lambda}$
of integers converging to $n_0$.  Define the {\it smooth functions on
$\beta({\Bbb Z})$}, denoted 
by $\ell^{\infty \infty}({\Bbb Z})$ or $C^\infty(\beta({\Bbb Z}))$, 
to be those functions $\varphi$ in $\ell^\infty({\Bbb Z})$ that 
satisfy
$$\lim_{\alpha \in \Lambda} 
n_\alpha^d \biggl(\varphi(n_\alpha) - \varphi(n_0)\biggr) = 0
\eqno {(2.2)} $$
for each $d=0, 1, 2, \dots$ and for each $n_0 \in \beta({\Bbb Z})$.
This set of functions $\ell^{\infty \infty}({\Bbb Z})$ 
is a dense $\star$-subalgebra of $\ell^\infty({\Bbb Z})$,
which contains $\ell^\infty_c({\Bbb Z})$, 
and plays an analogous role to 
$\ell^\infty({\Bbb Z})$ as
${\cal S}({\Bbb Z})$, the Schwartz 
functions on $\Bbb Z$, does to $c_0({\Bbb Z})$.
If $p$ is a projection 
in $\ell^\infty({\Bbb Z})$, we noticed above that $p(n_0)=0$
or $p(n_0)=1$ for any 
$n_0 \in \beta({\Bbb Z})$.  In either case, the quantity in
parentheses in (2.2) eventually becomes $0$, so that $\varphi = p$
satisfies (2.2).  This 
shows that $\ell^\infty_c({\Bbb Z}) \subseteq \ell^{\infty \infty}({\Bbb Z})$.

\noindent{\bf Lemma 2.3. \ }  The limit (2.2) holds independently of the
choice of net $\{ n_\alpha \}_{\alpha \in \Lambda}$ converging to
$n_0$.

\noindent{\bf Proof:} Let $\epsilon > 0 $ be given and find a $\beta$
such that 
$ \bigl| n_\alpha^d \bigl( \psi(n_\alpha) - \psi(n_0) \bigr) \bigr| 
< \epsilon $
for $\alpha \ge \beta$.  The set 
$S=\bigcup \{ n_\alpha \mid \alpha \ge \beta \}$
is in the ultrafilter associated with $n_0$ and 
$ \bigl| m^d \bigl( \psi(m) - \psi(n_0) \bigr) \bigr| 
< \epsilon $
for $m \in S$.  If $\{ m_\alpha \}_{\alpha \in \Gamma}$ is another net
tending to $n_0$, then it is eventually in $S$.  So we have
$ \bigl| m_\alpha^d \bigl( \psi(m_\alpha) - \psi(n_0) \bigr) \bigr| 
< \epsilon $
for $\alpha \ge \gamma$, for some $\gamma \in \Gamma$.  {\bf QED}

To see that $\ell^{\infty \infty}({\Bbb Z})$ is closed under products, let $\varphi$, $\psi \in \ell^{\infty \infty}({\Bbb Z})$.
Then evaluate the quantity in parentheses in (2.2), namely the difference
$$\varphi(n_\alpha)  \psi(n_\alpha)  - \varphi(n_0)  \psi(n_0) = 
	 \biggl( \varphi(n_\alpha) - \varphi(n_0) \biggr) \psi(n_\alpha)
	+ \varphi(n_0)  \biggl( \psi(n_\alpha)  - \psi(n_0) \biggr).$$
So the absolute value of the quantity in the limit (2.2) is 
$$\biggl| n_\alpha^d \biggl( 
		\varphi\psi(n_\alpha)-\varphi\psi(n_0)\biggr)\biggr|\le
 \biggl| n_\alpha^d \biggl( \varphi(n_\alpha) - \varphi(n_0) \biggr) \biggr|
 \| \psi \|_\infty
	+ \| \varphi \|_\infty  
 \biggl| n_\alpha^d \biggl( \psi(n_\alpha) - \psi(n_0) \biggr) \biggr|
$$
Clearly, this tends to zero as $n_\alpha$ tends to $n_0$.

Next, we note that $\ell^{\infty \infty}({\Bbb Z})$ is actually bigger than 
$\ell^\infty_c({\Bbb Z})$.
For example, any function in ${\cal S}({\Bbb Z})$ satisfies the limit (2.2), so
$\ell^{\infty \infty}({\Bbb Z}) 
\supseteq \ell^\infty_c({\Bbb Z}) + {\cal S}({\Bbb Z})$.  
For $\varphi \in {\cal S}({\Bbb Z})$, note that
$\varphi(n_0) = 0$ for any $n_0 \in \beta({\Bbb Z}) - {\Bbb Z}$.  
(Any nonprincipal ultrafilter eventually leaves every finite set.)
Therefore
$$\biggl| n_\alpha^d \biggl( 
		\varphi(n_\alpha)-\varphi(n_0)\biggr)\biggr| = 
 \biggl| n_\alpha^d  \varphi(n_\alpha) \biggr| = 
 1/n_\alpha^2 \biggl| n_\alpha^{d+2}  \varphi(n_\alpha) \biggr| \le
 1/n_\alpha^2 \|  \varphi \|_{d+2}, 
$$
where $\| \quad \|_{d+2}$ 
denotes the $d+2$th Schwartz seminorm on ${\cal S}({\Bbb Z})$.

\noindent{\bf Proposition 2.4. \ }  The inclusions 
$\ell^\infty_c({\Bbb Z}) + {\cal S}({\Bbb Z}) 
\subseteq \ell^{\infty \infty}({\Bbb Z}) 
\subseteq \ell^\infty({\Bbb Z})$ are proper.

\noindent{\bf Proof:} 
The function 
$1\over{n^2 + 1}$ in $c_0({\Bbb Z}) - {\cal S}({\Bbb Z})$ 
does not satisfy (2.2),  
so the second inclusion is proper.
The first inclusion is proper 
since $e^{i{\cal S}({\Bbb Z})} \subseteq \ell^{\infty \infty}({\Bbb Z})$.
{\bf QED}

Let $\varphi \in \ell^{\infty \infty}({\Bbb Z})$.  We may write $\varphi$ (in fact any function
in 
$\ell^\infty({\Bbb Z})$) uniquely in the form
$$
\varphi = \sum_{q=1}^\infty c_q p_q, \eqno (2.5) 
$$
where the $p_q$'s are {\it pairwise disjoint projections}, and
the $c_q$'s are {\it distinct} constants.   

\noindent{\bf Proposition 2.6. \ }  At most finitely many projections in the
series (2.5) have infinite support. 

\noindent{\bf Proof:} 
{\it We assume for a contradiction that
there are infinitely many projections in (2.5), each having 
infinite support}.  
%By replacing $\varphi$ with its absolute value, 
%and then scaling down by the constant $\pa \varphi \pa_\infty$,
%we may assume that $\varphi$ has unit norm and is non-negative.
%Then $0 \le c_q \le 1$ for each $q$.  (When we replaced $\varphi$ 
%with its absolute value, we may have had to combine some of the projections 
%$p_q$ into one  larger projection, in order to keep the $c_q$'s distinct.
%There will still be infinitely many projections
%with infinite support after we do this.)
Note that the 
coefficients of the infinite projections must have an accumulation
point.  By dropping to a subsequence, and getting rid of all the 
finite projections, we may assume that $c_q \rightarrow c_0 \in {\Bbb C}$ 
as $q \rightarrow \infty$.  (Multiply $\varphi$ by an appropriate projection,
and renumber the $c_q$'s.)
Let $S_q\subseteq {\Bbb Z}$ denote the support of the projection $p_q$.
Then each set $S_q$ is infinite.
Define a decreasing sequence of infinite subsets of $\Bbb Z$ by the 
disjoint unions
$$U_n = \bigcup_{q \ge n} \biggl( S_q  \bigcap 
\biggl\{ m \in {\Bbb Z}\,\, \biggl| \,\, m  
\ge {1\over { | c_q - c_0 | }} \biggr\}  \biggr)  $$
for $n = 1, 2, \dots $.   
Since finite intersections of the sets $U_n$'s are non-empty, there exists
an ultrafilter $\cal U$ for which $U_n \in {\cal U}$ for each $n$.
Let $\{ n_\alpha \}_{\alpha \in \cal U}$ be a net of integers that is 
eventually in this ultrafilter.  This net must therefore also eventually be in 
each of the sets $U_n$.  By construction, the point $n_0 \in {\cal FU}({\cal U})
\in \beta({\Bbb Z})$ that $\{ n_\alpha \}_{\alpha \in {\cal U}}$ converges to 
must satisfy $\varphi(n_0)= c_0$.

If $n_\alpha \in S_q \bigcap \{ m \in {\Bbb Z} \mid m \ge {1\over |c_q - c_0|}
\}$, then
\begin{eqnarray}
\biggl| n_\alpha^d \biggl( \psi(n_\alpha) - \psi(n_0) \biggr) \biggr| & =
\biggl| n_\alpha^d \biggl( c_q - c_0 \biggr) \biggr| 
\qquad \qquad {\rm since} \quad  n_\alpha \in S_q
\nonumber
\\
& \ge \bigl| n_\alpha^d \times {1\over {n_\alpha}} \bigr| 
= \bigl| n_\alpha^{d-1} \bigr| 
\ge 1, 
\nonumber
\end{eqnarray}
contradicting the fact that $\varphi$ must satisfy (2.2) for $d \ge 2$
at the point $n_0 \in \beta({\Bbb Z})$ we constructed.  {\bf QED}

It follows from Proposition 2.6 that for any $\varphi \in \ell^{\infty \infty}({\Bbb Z})$,
we may substract an element of $\ell^\infty_c({\Bbb Z})$ 
to force the expansion (2.5)
to have only projections of finite support.

\noindent{\bf Remark 2.7. \ }  Note that ${\cal S}({\Bbb Z})$ is an ideal in $\ell^{\infty \infty}({\Bbb Z})$ and 
$\ell^\infty({\Bbb Z})$.  The closure 
${\overline {{\cal S}({\Bbb Z})}}^{\| \,\,\,\, \|_\infty}$
is equal to $c_0({\Bbb Z})$,
so ${\cal S}({\Bbb Z})$ is not dense in $\ell^\infty({\Bbb Z})$.  
The algebra $\ell^\infty_c({\Bbb Z})$, being unital,
is {\it not} an ideal in 
either algebra $\ell^{\infty \infty}({\Bbb Z})$ 
or $\ell^\infty({\Bbb Z})$, and for the same
reason 
$\ell^{\infty \infty}({\Bbb Z})$ is 
not an ideal in $\ell^\infty({\Bbb Z})$.  This is in contrast to the
case $c_c({\Bbb Z}) \subseteq {\cal S}({\Bbb Z}) \subseteq c_0({\Bbb Z})$, 
where every algebra 
is a dense ideal in every algebra above it.

\vskip\baselineskip
\vskip\baselineskip
\section{References}
\smallskip
\smallskip

\noindent[{\bf Roy, 1968}]
\, H.L. Royden,
{\it Real Analysis}, second edition,
Macmillan Publishing Co., Inc., New York, 1968.

\noindent[{\bf Sch, 1998}]
\, L. B. Schweitzer,
{\it $C^{\infty}$ functions on the cantor set, and a smooth $m$-convex
Fr\'echet subalgebra of $O_2$},
Pacific J. Math. {\bf 184(2)} (1998), 349-365.

\vskip\baselineskip
\noindent{Email: lsch@svpal.org.
Web Page: http://www.svpal.org/$\thicksim$lsch/Math.}

\end{document}